\documentclass[11pt, reqno]{amsart}
\usepackage{amssymb, amsmath, amsthm, amsfonts}

\usepackage{graphics}

\newcommand{\demo}{\begin{proof}}
\newcommand{\edemo}{\end{proof}}
\newcommand{\demoname}[1]{\begin{proof}[#1]}
\newcommand{\edemoname}{\end{proof}}

\theoremstyle{plain}
\newtheorem{theorem}{Theorem}[section]
\newtheorem{conjecture}[theorem]{Conjecture}
\newtheorem{corollary}[theorem]{Corollary}
\newtheorem{lemma}[theorem]{Lemma}

\theoremstyle{definition}
\newtheorem{example}[theorem]{Example}
\newtheorem*{definition}{Definition}
\newtheorem{remark}[theorem]{Remark}

\newcommand{\thm}{\begin{theorem}}
\newcommand{\ethm}{\end{theorem}}
\newcommand{\conj}{\begin{conjecture}}
\newcommand{\econj}{\end{conjecture}}

\newcommand{\expl}{\begin{example}}
\newcommand{\eexpl}{\qex\end{example}}
\newcommand{\defn}{\begin{definition}}
\newcommand{\edefn}{\qef\end{definition}}
\newcommand{\remk}{\begin{remark}}
\newcommand{\eremk}{\qex\end{remark}}

\newcommand{\coro}{\begin{corollary}}
\newcommand{\ecoro}{\end{corollary}}
\newcommand{\lem}{\begin{lemma}}
\newcommand{\elem}{\end{lemma}}

\providecommand{\qexsymbol}{$\lozenge$}%
\newcommand{\mathqex}{\quad\hbox{\qexsymbol}}
\DeclareRobustCommand{\qex}{%
  \ifmmode \mathqex
  \else
    \leavevmode\unskip\penalty9999 \hbox{}\nobreak\hfill
    \quad\hbox{\qexsymbol}%
  \fi
}
\providecommand{\qefsymbol}{$\triangle$}%
\newcommand{\mathqef}{\quad\hbox{\qefsymbol}}
\DeclareRobustCommand{\qef}{%
  \ifmmode \mathqef
  \else
    \leavevmode\unskip\penalty9999 \hbox{}\nobreak\hfill
    \quad\hbox{\qefsymbol}%
  \fi
}

\newcommand{\enum}{\begin{enumerate}}
\newcommand{\eenum}{\end{enumerate}}
\newcommand{\cset}[2]{C^{{#1}}_{{#2}}}
\newcommand{\spce}{\smallskip}
\newcommand{\rrr}[1]{{\mathbb{R}}^{#1}}

\newcommand{\nd}[1]{$#1$\nobreakdash-\hspace{0pt}}
\newcommand{\kmnnp}[2]{D_{{#1}}({#2})}
\newcommand{\kmnna}[2]{\tilde D_{{#1}}({#2})}
\newcommand{\ka}[2]{D_{{#1}}^a({#2})}
\newcommand{\stcurv}{standard curvature}
\newcommand{\sdcurv}{stratified curvature}

\newcommand{\sdcurva}{ascending stratified curvature}

\newcommand{\modsdcurv}{modified stratified curvature}

\newcommand{\sddef}{generalized angle defect}
\newcommand{\sddefs}{generalized angle defects}
\newcommand{\Sddef}{Generalized angle defect}

\newcommand{\seccc}{stratified Euler characteristic}
\newcommand{\aang}{\alpha}

\newcommand{\eucn}{\chi^s}
\newcommand{\simpn}[1]{\nd{{#1}}dimensional simplicial complex}

\newcommand{\faceof}{\prec}
\newcommand{\offace}{\succ}

\newcommand{\simrankk}[2]{T_{{#1}}({#2})}
\newcommand{\strank}{stratification rank}

\newcommand{\ads}{angle defect sequence}
\newcommand{\pnr}[2]{P_{{#1}, {#2}}}

\DeclareMathOperator{\starpl}{star}
\DeclareMathOperator{\linkpl}{link}

\begin{document}

\title{The angle defect for odd-dimensional simplicial manifolds}
\author{Ethan D.\ Bloch}
\address{Bard College\\
Annandale-on-Hudson, NY 12504\\
U.S.A.}
\email{bloch@bard.edu}
\thanks{The author would like to thank Mark Halsey, Lauren Rose and Rebecca Thomas for their assistance}
\date{}
%%[[52B70 is polyhedral manifolds]]
\begin{abstract}
In a 1967 paper, Banchoff stated that a certain type of polyhedral curvature, that applies to all finite polyhedra, was zero at all vertices of an odd-dimensional polyhedral manifold; one then obtains an elementary proof that odd-dimensional manifolds have zero Euler characteristic.  In a previous paper, the author defined a different approach to curvature for arbitrary simplicial complexes, based upon a direct generalization of the angle defect.  The generalized angle defect is not zero at the simplices of every odd-dimensional manifold.  In this paper we use a sequence based upon the Bernoulli numbers to define a variant of the angle defect for finite simplicial complexes that still satisfies a Gauss-Bonnet type theorem, but is also zero at any simplex of an odd-dimensional simplicial complex $K$ (of dimension at least $3$), such that $\chi(\linkpl(\eta^i, K)) = 2$ for all \nd{i}simplices $\eta^i$ of $K$, where $i$ is an even integer such that $0 \le i \le n-1$.  As a corollary, an elementary proof is given that any such simplicial complex has Euler characteristic zero.
\end{abstract}
\maketitle

\section{Introduction}
\label{secINT}

For a triangulated polyhedral surface $M^2$, the usual 
notion of curvature at a vertex $v$ is the classical angle defect $d_v 
= 2\pi - \sum \alpha_i$, where the $\alpha_i$ are the angles 
of the triangles containing $v$.  This curvature function, 
which goes back at least as far as Descartes (see\cite{FE}), 
satisfies all the standard properties one would expect a 
curvature function on polyhedra to satisfy.  For example, 
the angle defect is invariant under local polyhedral 
isometries (that is, functions that preserve the lengths of 
edges); it is locally defined; it is zero at a vertex that has a flat star; it is invariant under subdivision; and it 
satisfies the polyhedral Gauss-Bonnet Theorem, which says 
$\sum_v d_v = 2\pi \chi(M^2)$, where the summation is over 
all the vertices of $M^2$, and $\chi(M^2)$ is the Euler 
characteristic of $M^2$.

There are two approaches to generalizing the classical angle defect to arbitrary (finite) polyhedra in all dimensions.  One method, which we will refer to as \stcurv, has been studied extensively from a differential geometric point of view.  See, among others, \cite{BA1}, \cite{WINT}, \cite{BUDA}, \cite{CHEE}, \cite{C-M-S} and \cite{ZA}.  This approach to generalizing the angle defect, which is based on exterior angles, is simple to define, and it's convergence properties has been well studied.  In \stcurv, all the curvature is concentrated at the vertices, as in the case for the classical angle defect of polyhedral surfaces.  Another approach, called simply the angle defect (also known as the angle deficiency), has been studied in the case of convex polytopes by a number of combinatorialists, for example \cite{SH2} and \cite{GR2}; more generally, for the wider study of angle sums in convex polytopes, see for example \cite[Chapter~14]{GR1}, \cite{SH}, \cite{P-S} and \cite{MC}.  In \cite{G-S2} a Gauss-Bonnet type theorem (also referred to as Descartes' Theorem) is proved for the angle defect in polytopes with underlying spaces that are manifolds.  In contrast to \stcurv, which is concentrated at vertices, the angle defect for convex polytopes is found at each simplex of codimension at least $2$ (it can be defined for all simplices, but the angle defect at a codimension $0$ or $1$ simplex will always be zero).   

In \cite{BL4} we extended the notion of angle defect to arbitrary simplicial complexes, not just manifolds, by using a simple toplogical decomposition of each simplicial complex.  In order to compare our approach with \stcurv, we originally concentrated our curvature at the vertices, and called it \sdcurv\; see \cite[Section~3]{BL4} for details.  In \cite[Section~4]{BL4}, we took the variant approach most directly comparable to the combinatorial authors listed above, in that we left the pure angle defects defined for each simplex of codimension at least $2$.  In \cite[Section~4]{BL4} we referred to this approach by the unfortunate name ``\modsdcurv,'' which really misses the point that in this approach we are really still working with a pure angle defect.  Hence, in the present paper, we will use the better name ``\sddef'' (which we also use in \cite{BL10}).   

A detailed comparison of \stcurv\ with both \sdcurv\ and the \sddef\ may be found in \cite[Section~4]{BL4}.  We mention here, however, that all these approaches satisfy some of the basic properties that one would expect of curvature, such as being locally defined, invariant under local isometries, and satisfying a Gauss-Bonnet type theorem (though the Gauss-Bonnet Theorem for \sdcurv\ and the \sddef\ uses a modified Euler characteristic rather than the standard Euler characteristic, as discussed in \cite[Section~2]{BL4}).  In \cite{BL10} we show that the \sddef\ has a Morse theoretic interpretation very similar to (though not quite as simple as) the Morse theoretic approach to \stcurv\ found in \cite{BA1}, \cite{BA2} and \cite{BA3}.  Hence, on most counts, it is fair to say that the \sddef\ behaves as nicely as \stcurv.  

There is one place, however, where the \sddef\ falls short of \stcurv.  In \cite[Section~5]{BA1}, Banchoff indicates that for an odd-dimensional polyhedral manifold, the \stcurv\ is zero at every vertex (the proof of the main step of this claim is not given, however).  It follows from the Gauss-Bonnet theorem for \stcurv\ that every odd-dimensional polyhedral manifold has Euler characteristic zero (a well-known fact, but Banchoff's method yields a completely elementary proof of this fact, without using algebraic topology).  
By contrast, as we will see in Section~\ref{secSTA}, there are odd-dimensional compact simplicial manifolds for which neither \sdcurv\ nor the \sddef\ is identically zero.

The purpose of the present paper is to define a variant on the \sddef\ called \sdcurva\, such that \sdcurva\ satisifies most of the nice properties of the \sddef, and yet also has the additional property that it is identically zero for any odd-dimensional simplicial complex $K$ (of dimension at least $3$), such that $\chi(\linkpl(\eta^i, K)) = 2$ for all \nd{i}simplices $\eta^i$ of $K$, where $i$ is an even integer such that $0 \le i \le n-1$.  We then deduce that any such simplicial complex has Euler characteristic zero.  In particular, these results hold for any odd-dimensional compact simplicial homology manifold.

The outline of the paper is as follows.  For the sake of completeness, we start off in Section~\ref{secSSC} with a very brief review of the needed definitions and theorems from \cite{BL4}, leaving all the details to that paper.  We give our new definitions and statements of results in Section~\ref{secSTA}, and then give proofs in Section~\ref{secPRO}.

Throughout this paper we will restrict our attention to simplicial complexes, rather than more general polyhedra.  We make the following assumptions, which we use throughout this paper without stating: All simplicial complexes are finite, and are embedded in Euclidean space (which we will not name when it is not necessary); all manifolds are compact and without boundary.

\section{Review of the \Sddef}
\label{secSSC}

We give here a very brief review of those definitions and statements of results from \cite{BL4} that we will be using, leaving proofs and examples to the original paper.

Throughout this section, let $K$ be an \simpn n.  Whereas in \cite{BL4} we allow for a certain class of non-embedded simplicial complexes, here for convenience we look only at actual simplicial complexes in Euclidean space.  If $\eta$ and $\sigma$ are simplices in $K$, we write $\eta 
\faceof \sigma$ to indicate that $\eta$ is a face of $\sigma$.  As usual, we let $|K|$ denote the underlying 
space of $K$.

For convenience (though not necessity), we adopt the convention that all angles are normalized so that the volume of the unit \nd{(n-1)}sphere in \nd{(n-1)}measure is $1$ in all dimensions.  For any \nd{n}simplex $\sigma^n$ in Euclidean space, and any \nd{i}face $\eta^i$ of $\sigma^n$, let $\aang(\eta^i, \sigma^n)$ denote the solid angle in $\sigma^n$ along $\eta^i$, where by normalization such an angle is always a number in the interval $[0, 1]$.  

\remk\label{remANG1}
Let $\sigma^n$ be an \nd{n}simplex.  Observe that $\aang(\eta^{n-1}, \sigma^n) = 1/2$ for any \nd{(n-1)}face $\eta^{n-1}$ of $\sigma^n$, and that $\aang(\sigma^n, \sigma^n) = 1$.
\eremk

\begin{figure}[ht]
\hfil\includegraphics{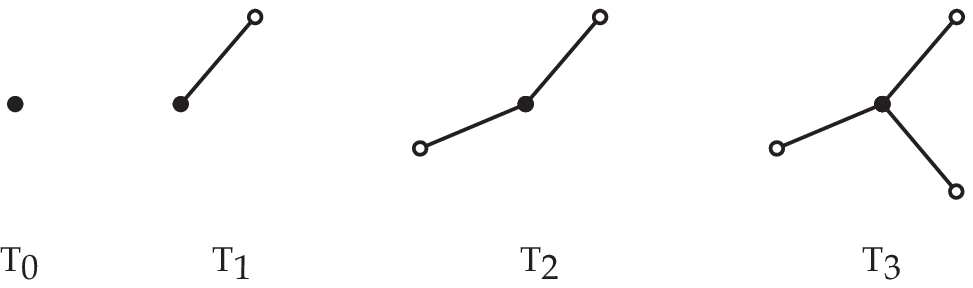}\hfil
\caption{}\label{figPH21}
\end{figure}

\begin{figure}[ht]
\hfil\includegraphics{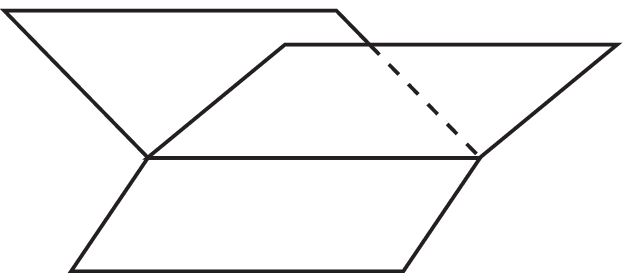}\hfil
\caption{}\label{figPH22}
\end{figure}

\defn 
For each non-negative integer $i$, let $T_i$ denote the open cone on $i$ points; alternately, $T_i$ is the space obtained by gluing together $i$ copies of the half open interval $[0, 1)$ at the point $\{ 0\}$ in each.  We take $T_0$ to be a single point.  See Figure~\ref{figPH21}.  Let $\pnr ni$ denote the space $\pnr ni = T_i \times \rrr{n-1}$.  If $*$ denotes the cone point of $T_i$, we call $\{*\} \times \rrr{n-1} \subseteq \pnr ni$ the {\bf apex set} of $T_i$.  See Figure~\ref{figPH22}.
\edefn

Observe that $\pnr ni$ is not homeomorphic to $\pnr nj$ when $i \neq j$.

We now need to think of simplices as open 
(and hence disjoint).  

\defn 
Let $K$ be an \simpn n.  For each 
non-negative integer $r$ such that $r \ne 2$, we define the 
subset $\cset nr(K)$ of $|K|$ by 
\begin{align*}
\cset nr(K) &= \{ x \in |K| \mid x \text{ has nbhd.\
homeomorphic to }\pnr nr, \\
&\qquad \text{ where the homeomorphism takes $x$ to the apex set of } \pnr nr\}.
\end{align*}
Define
$$
\cset n2(K) = |K| - \bigcup_{r \ne 2} \cset nr(K).
$$
\edefn

\remk\label{remCSET}
(1) Let $K$ be an \simpn n.  The sets $\cset nr(K)$ are well 
defined, since each $x \in |K|$ can have a neighborhood 
homeomorphic to $\pnr nr$ (where the homeomorphism takes $x$ 
to the apex set of $\pnr nr$) for at most one number $r$.  Moreover, the sets $\cset nr(K)$ are well 
defined up to homeomorphism of $|K|$.
\spce

\noindent (2) Because $K$ is a finite simplicial complex, 
there is some positive integer $P$ such that $\cset nr(K) = 
\emptyset$ for all $r > P$.
\spce

\noindent (3) The sets $\cset nr(K)$ are disjoint, and cover 
$|K|$.  For each $r \neq 2$, the set $\cset nr(K)$ is an 
$(n-1)$-manifold.  Moreover, each set $\cset nr(K)$ is the 
union of (open) simplices of $K$, since all points in any 
simplex of $K$ have homeomorphic neighborhoods in $|K|$ (if 
the neighborhoods are taken small enough).  If $\sigma \in 
K$, then $\sigma \subseteq \cset nr(K)$ for some unique integer 
$r$.
\eremk

\defn 
Let $K$ be an \simpn n.  Let $\sigma 
\in K$ be a simplex.  The {\bf \strank} of $\sigma$, denoted 
$\simrankk n\sigma$, is defined by $\simrankk n\sigma = r/2$, 
where $\sigma \in \cset nr(K)$ for some unique integer $r$.
\edefn

\remk\label{remCURV2}
Let $K$ be an \simpn n.  If $\eta^{n-1}$ is an \nd{(n-1)}simplex of $K$, then $\simrankk n{\eta^{n-1}}$ equals one half the number of \nd{n}simplices of $K$ that have $\eta^{n-1}$ as a face.  If $\sigma^n$ is an \nd{n}simplex of $K$, then $\simrankk n{\sigma^n} = 1$.
\eremk

The following definition was originally given in \cite[Section~4]{BL4}, though here we use the better name given given below (and also used in \cite{BL10}, as discussed in Section~\ref{secINT}).  As in \cite{SH2} and \cite{GR2}, for example, the curvature as we discuss it is to be found at all simplices (though the non-zero curvature is always at simplices of codimension at least $2$), in contrast to the approach of \cite{BA1} et al., where the curvature is concentrated at the vertices.

\defn
Let $K$ be an \simpn n, and let $\eta^i$ be an \nd{i}simplex of $K$, where $0 \le i \le n$.  The {\bf \sddef} at $\eta^i$ is the number $\kmnna n{\eta^i}$ defined by
\begin{equation}\label{eqMOD1}
\kmnna n{\eta^i} = \simrankk n{\eta^i} - \sum_{\sigma^n \offace \eta^i} \aang(\eta^i, \sigma^n),\end{equation}
where the summation is over all \nd{n}simplices $\sigma^n$ which have $\eta^i$ as a face.
\edefn

\remk\label{remCURV1}
Using Equation~\ref{eqMOD1}, together with Remarks~\ref{remANG1} and \ref{remCURV2}, it is seen that $\kmnna n{\eta^i} = 0$ if $i = n-1$ or $i = n$.
\eremk

In \cite[Section~3]{BL4} we gave the following variant definition of curvature of simplicial complexes; this definition concentrates all curvature at the vertices, which is somewhat unnatural, but was used to facilitate a comparison with the definition of curvature found in \cite{BA1} et al.  

\defn
Let $K$ be an \simpn n, and let $v$ be a vertex of $K$.  The {\bf \sdcurv} at $v$ is the number $\kmnna n{\eta^i}$ defined by
\begin{equation}\label{eqMOD2}
\kmnnp n{v, K} = \sum_{i=0}^{n-2} \frac{(-1)^i}{i+1} \sum_{\eta^i \ni v} \kmnna n{\eta^i}.
\end{equation}
\edefn

In \cite{BL4} we proved Gauss-Bonnet type theorems for both \sddef\ and \sdcurv, though rather than using the standard Euler characteristic, we used the following variant of the Euler characteristic.  We will use this new characteristic in the present paper as well.

\defn
Let $K$ be an \simpn n.  The {\bf \seccc} of $K$, denoted $\eucn (K)$, is defined by
\begin{equation}\label{eqSECCC}
\eucn (K) = \sum_{\eta \in K} \simrankk n\eta \,(-1)^{\dim \eta}.
\end{equation}
\edefn

\remk\label{remCURV4}
As discussed in \cite[Section~2]{BL4}, the \seccc\ of a simplicial complex is a homeomorphism invariant of its underlying space, though not a homotopy type invariant.
\eremk

\section{Statement of Results}
\label{secSTA}

We start by showing that  there are odd-dimensional compact simplicial manifolds for which neither \sdcurv\ nor the \sddef\ is identically zero.  We will be making use of convex polytopes in our discussion.  For general information on convex polytopes, see \cite{GR1}.  
As usual, if $Q$ is a simplicial complex, we let $f_i(Q)$ denote the number of \nd{i}simplices of $Q$.

To see that the \sddef\ is not identically zero for odd-dimensional simplicial manifolds is simple.  Let $K$ be the boundary of any \nd{n}dimensional convex polytope in $\rrr{n}$ with $n$ an even integer such that $n \ge 4$.  Hence $K$ is an \nd{(n-1)}dimensional simplicial sphere.  Let $\eta^i$ be an \nd{i}simplex of $K$ with $0 \le i \le n-2$.  Then by Theorem~(3) of \cite{SH2}, we know that $\kmnna 7{\eta^i} > 0$.  Hence the \sddef\ is certainly not zero for all simplices of all odd-dimensional combinatorial manifolds.  

We now turn to \sdcurv, where we need a slightly more complicated counterexample.  (The simplest counterexample we found is \nd{7}dimensional; it can be verified that there cannot be a \nd{3}dimensional counterexample, but the author does not know whether a \nd{5}dimensional counterexample could exist, though one cannot be constructed by our method.)

\expl\label{explSEVEN}
Take a triangle in the plane, and let $Q$ be the cone of the cone of the suspension of the triangle .  Then $Q$ is a \nd{5}dimensional convex polytope in $\rrr{5}$.  Let $L$ be the boundary of $Q$, so that $L$ is a \nd{4}dimensional simplicial sphere.  It is straightforward to see that $f_0(L) = 7$, $f_1(L) = 20$, $f_2(L) = 29$, $f_3(L) = 22$, and $f_4(L) = 8$.  Triangulate $Q$ by putting a vertex $v$ in the interior of $Q$, and then joining $v$ to $L$.

Next, let $M$ be the boundary of an \nd{8}dimensional simplicial convex polytope in $\rrr{8}$, so that $M$ is a \nd{7}dimensional simplicial sphere; it does not matter what simplicial convex polytope is used.  Let $\eta^5$ be some \nd{5}simplex in $M$.  As mentioned above, we know that $\kmnna 7{\eta^5} > 0$.  We subdivide $M$ as follows.  First, put a copy of the triangulation of $Q$ in the interior of $\eta^5$.  Then add further simplices to $M$, without subdividing $Q$, so that we obtain a simplicial complex, denoted $M'$.

Suppose that $\alpha^i$ is an \nd{i}simplex of $M'$.  Then $\alpha^i$ is contained in a unique \nd{r}simplex $\zeta^r$ of $M$, where $i \le r \le 7$ (recall that we are thinking of simplices as open, and hence disjoint).  It is seen that $\simrankk 7{\alpha^i} = \simrankk 7{\zeta^r}$ and $\kmnna 7{\alpha^i} = \kmnna 7{\zeta^r}$, where the left hand side of each of these equations is with respect to $M'$, and the right hand side is with respect to $M$.

Let $\omega^i$ be a simplex in $M'$ that has $v$ as a vertex.  If $\omega^i$ is contained in $\eta^5$ (and hence it is contained in $Q$), then $\kmnna 7{\omega^i} = \kmnna 7{\eta^5}$.  If $\omega^i$ is not contained in $\eta^5$, then it is contained in one of the \nd{6}dimensional or \nd{7}dimensional simplices of $M$ that have $\eta^5$ as a face.  By Remark~\ref{remCURV1}, it follows that $\kmnna 7{\omega^i} = 0$.  Now, using Equation~\ref{eqMOD2}, and keeping track of the simplices of $Q$ that contain $v$ by looking at $\linkpl(v, Q) = L$, we see that
\begin{align*}
\kmnnp 7{v, M'} &= \kmnna 7{\eta^5} - \frac 12f_0(L)\kmnna 7{\eta^5} + \frac 13f_1(L)\kmnna 7{\eta^5} - \frac 14f_2(L)\kmnna 7{\eta^5}\\
 &\qquad\qquad + \frac 15f_3(L)\kmnna 7{\eta^5} - \frac 16f_4(L)\kmnna 7{\eta^5}\\
 &= - \frac 1{60}\kmnna 7{\eta^5} \ne 0.
\end{align*}
Hence we have an example of an odd-dimensional simplicial manifold with non-zero \stcurv\ at one of its vertices.  (We note that this example is a combinatorial manifold, using the fact that all simplicial convex polytopes, and all subdivisions of simplicial convex polytopes, are combinatorial manifolds; this follows from standard results about combinatorial manifolds, as found in \cite[Section I.5]{HU}, for example). 
\eexpl

Our goal is to reformulate the notion of curvature for simplicial complexes in terms of the \sddef, with the aim of remedying the above problem; that is, we want a notion of curvature based on the \sddef, but which is identically zero at every simplex of odd-dimensional simplicial manifolds, and also satisfies the nice properties of the \sddef, such as a Gauss-Bonnet type theorem.  Certainly, the most natural way to look at the \sddef\ is simply to use $\kmnna n{\tau^p}$ as defined for each simplex $\tau^p$ in a simplicial complex.  In \sdcurv, we concentrated the $\kmnna n{\tau^p}$ at the vertices via Equation~\ref{eqMOD2}, to obtain $\kmnnp n{v, K}$.  Concentrating the curvature of a simplicial complex at the vertices isn't useful for our present purposes, but it raises the question of whether there is some other way of redistributing the \sddefs\ in such a way that the other desired properties of simplicial curvature still hold (for example, locally defined, invariant under local isometries, and a Gauss-Bonnet type theorem), and yet the property that the curvature is zero at all simplices in odd-dimensional simplicial manifolds also holds.  

What does it mean for curvature to be locally defined in a simplicial complex?  In the smooth case, for comparison, the fact that curvature is locally defined means that it can be calculated on an arbitrarily small neighborhood of each point.  In the simplicial case, a number associated to each simplex is local if it can be defined using only the star of each simplex.  In particular, let $K$ be an \simpn n, and let $\tau^p$ be a \nd{p}simplex of $K$.  Anything that can be computed in $\starpl (\tau^p, K)$ can be considered to be locally defined at $\tau^p$.  For example, clearly $\kmnna n{\tau^p}$ is locally defined at $\tau^p$.  However, if $\eta^i$ has $\tau^p$ as a face (so that $\eta^i \in \starpl (\tau^p, K)$), then $\starpl (\eta^i, K)$ is contained in $\starpl (\tau^p, K)$, and hence $\kmnna n{\eta^i}$ is also locally defined at $\tau^p$.  Therefore, we see that any number produced out of $\kmnna n{\tau^p}$ and all the $\kmnna n{\eta^i}$ for $\eta^i$ that have $\tau^p$ as a face is locally defined at $\tau^p$.  Of course, a number produced out of an arbitrarily chosen combination of \sddefs\ will not necessarily behave nicely as one would expect of a curvature function.  There is, we will see, one method of combining these \sddefs\ that has the desired properties.  

To define the desired curvature, we will need the following definition, which uses the Bernoulli numbers.  (See, for example, \cite[Sections~6.5 and 7.6]{G-K-P} for more on the Bernoulli numbers.)  As usual, we let $B_k$ denotes the $k^{\text{th}}$ Bernoulli number.  It is not entirely surprising to see the Bernoulli numbers appear here, given that they arise both in the study of angle sums of simplices (see \cite{PESC}, as well as \cite{GUIN}), and in topology (see \cite[Appendix B]{M-S}).

\defn
For each non-negative integer $n$, let $a_n$ be the number defined by
\begin{equation}\label{eqBERN1}
a_n = \frac {4B_{n+2}(2^{n+2} - 1)}{n+2}.
\end{equation}
We will refer to the sequence $a_0, a_1, a_2, \ldots$ as the {\bf \ads}.
\edefn

The numbers given in the above definition (modulo a factor of two) appear in the literature (for example \cite[p.\ 508]{KNOP}), though the sequence $a_0, a_1, a_2, \ldots$ does not appear to have a name, and hence we will use the name given above.  Like the Bernoulli numbers, the \ads\ does not appear to have an obvious pattern.  The first $18$ numbers in the \ads\ are
$$
1, 0, -\frac12, 0, 1, 0, -\frac{17}4, 0, 31, 0, -\frac{691}2, 0, 5461, 0, -\frac{929569}8, 0, 3202291, 0.
$$
The following lemma gives some properties of the \ads\ that we will use later on; the proof of the lemma will be given in Section~\ref{secPRO}.

\lem\label{lemADS1}
$\ $
\enum
\item[(1)] For every odd integers $n$, we have $a_n = 0$.
\item[(2)] The \ads\ satisfies the recursion relation $a_0 = 1$, and
\begin{equation}\label{eqADS1}
a_n + \sum_{i=0}^{n-1} \frac {a_i}2 \binom {n+1}{i+1} = 1
\end{equation}
for all positive integers $n$.
\eenum
\elem

We can now give the definition of our new type of curvature.

\defn
Let $a_0, a_1, a_2, \ldots$ be the \ads.  Let $K$ be an \simpn n, and let $\tau^p$ be a \nd{p}simplex of $K$.  The {\bf \sdcurva} at $\tau^p$ is the number $\ka n{\tau^p}$ defined by
\begin{equation}\label{eqCURVDEF1}
\ka n{\tau^p} = a_p\kmnna n{\tau^p} + \frac {a_p}2 \sum_{i=p+1}^{n} (-1)^{i-p} \sum_{\eta^i \offace \tau^p} \kmnna n{\eta^i}.
\end{equation}
\edefn

\remk\label{remCURV3}
(1) Using Remark~\ref{remCURV1}, it is seen that $\ka n{\tau^p} = 0$ if $p = n-1$ or $p = n$.  
Moreover, when $0 \le p \le n-2$, we have
\begin{equation}\label{eqCURVDEF2}
\ka n{\tau^p} = a_p\kmnna n{\tau^p} + \frac {a_p}2 \sum_{i=p+1}^{n-2} (-1)^{i-p} \sum_{\eta^i \offace \tau^p} \kmnna n{\eta^i}.
\end{equation}

\noindent (2) Because $a_p = 0$ for all odd $p$ (by Lemma~\ref{lemADS1} (1)), it follows that $\ka n{\tau^p} = 0$ for all odd $p$.  Hence, the \sdcurva\ is concentrated at the 
even-dimensional simplices.  (We could have defined \sdcurva\ only for even-dimensional simplices; however, it is more convenient to define the curvature are we did, including the odd-dimensional simplices.)
\eremk

It is straightforward to check that \sdcurva\ reduces to the classical angle defect when $K$ is a simplicial surface.

Our goal is to show that \sdcurva\ is nicely behaved.  Clearly it is locally defined, and invariant under local isometries.  Our three main theorems, which we now state, say that \sdcurva\ satisfies a Gauss-Bonnet type theorem (using the \seccc), that it  is identically zero for 
odd-dimensional compact simplicial manifolds, and that it is somewhat invariant under subdivision (to be clarified below).  The proofs of these results are given in Section~\ref{secPRO}.  Our first theorem is the following.

\thm\label{thmMAIN1}
Let $K$ be an \simpn n, where $n \ge 2$.  Then
\begin{equation}\label{eqMAIN1}
\sum_{\tau^p \in K} (-1)^p \ka n{\tau^p} = \eucn (K).
\end{equation}
\ethm

The following theorem gives a general condition that guarantees that \sdcurva\ will be identically zero in the odd-dimensional case.

\thm\label{thmMAIN2}
Let $K$ be an \simpn n, where $n$ is an odd integer such that $n \ge 3$.  Suppose that $\chi(\linkpl(\eta^i, K)) = 2$ for all \nd{i}simplices $\eta^i$ of $K$, where $i$ is an even integer such that $0 \le i \le n-1$.  Then $\ka n{\tau} = 0$ for every simplex $\tau$ of $K$.
\ethm

The following corollary is deduced immediately from Theorem~\ref{thmMAIN2}, using the fact that the links of simplices in homology manifolds have the appropriate homology groups (and hence Euler characteristic); see \cite[Theorem~63.2]{MU3} for more details.  Recall that we are assuming that all simplicial complexes are finite, and all manifolds are compact and without boundary.

\coro\label{coroMAIN3}
Let $K$ be a simplicial homology \nd{n}manifold, where $n$ is an odd integer such that $n \ge 3$.  Then $\ka n{\tau} = 0$  for every simplex $\tau$ of $K$.
\ecoro

The following result is another simple corollary to Theorem~\ref{thmMAIN2}.  

\coro\label{thmCOROMAIN2}
Let $K$ be an \simpn n, where $n$ is an odd integer such that $n \ge 3$.  Suppose that $\chi(\linkpl(\eta^i, K)) = 2$ for all \nd{i}simplices $\eta^i$ of $K$, where $i$ is an even integer such that $0 \le i \le n-1$.  Then $\chi(K) = 0$.
\ecoro

\demo
Combining Theorems~\ref{thmMAIN2} and \ref{thmMAIN1}, it follows that $\eucn (K) = 0$.  The condition on the links of even-dimensional simplices implies that every \nd{(n-1)}simplex of $K$ is the face of precisely two \nd{n}simplices.  It follows from Remark~\ref{remCURV2} that $\simrankk n{\tau^{p}} = 1$ for all \nd{p}simplices $\tau^{p}$ of $K$, where $p = n-1$ or $p = n$.  Moreover, by using Lemma~\ref{lemPSEU} (i), stated and proved below, it follows that $\simrankk n{\tau^{p}} = 1$ for all \nd{p}simplices $\tau^{p}$ of $K$, where $0 \le p \le n-2$.   It now follows from Equation~\ref{eqSECCC} that $\eucn(K) = \chi (K)$.
\edemo

Corollary~\ref{thmCOROMAIN2} immediately implies the following standard result.

\coro\label{thmCOROMAIN4}
Let $K$ be a triangulated homology \nd{n}manifold, where $n$ is an odd integer such that $n \ge 3$.  Then $\chi(K) = 0$.
\ecoro

Corollary~\ref{thmCOROMAIN4} is usually proved by algebraic topology, though our method yields a completely elementary proof of this fact (a similar type of proof is mentioned in \cite[Section~5]{BA1}).  However, it is not immediately evident how the standard proof using algebraic topology, which works for compact odd-dimensional manifolds, could be generalized to prove Corollary~\ref{thmCOROMAIN2}. 

The one remaining property of \sdcurva\ we wish to examine is invariance under subdivision.  Let $K$ be an \simpn n, where $n \ge 2$, and let $L$ be a subdivision of $K$.  Suppose that $\tau^s$ is an \nd{s}simplex of $L$.  Then $\tau^s$ is contained in a unique \nd{s}simplex $\zeta^p$ of $K$, where $s \le p \le n$ (recall that we are thinking of simplices as open, and hence disjoint).  To say that \sdcurva\ is invariant under subdivision would, in its nicest form, mean precisely that $\ka n{\tau^p} = \ka n{\zeta^r}$, where the left hand side of this equation is with respect to $L$, and the right hand side is with respect to $K$.  Because of the use of the \ads\ in the definition of \sdcurva, we do not quite have this simple equation, but we have a slightly modified version of it instead, as stated in the following theorem.

\thm\label{thmMAIN3A}
Let $K$ be an \simpn n, where $n \ge 2$, and let $L$ be a subdivision of $K$.  Let $\tau^s$ be an \nd{s}simplex of $L$ that is contained in a \nd{p}simplex $\zeta^p$ of $K$.  Then 
\begin{equation}\label{eqMAIN3A}
a_p\ka n{\tau^s} = a_s\ka n{\zeta^p}.
\end{equation}
When $s = p$, then 
\begin{equation}\label{eqMAIN3B}
\ka n{\tau^p} = \ka n{\zeta^p}.
\end{equation}
\ethm

For comparison, we note that both \stcurv\ and the \sddef\ are invariant under subdivision (as stated in \cite[Section~4]{BL4}), with no need for the modification introduced in Equation~\ref{eqMAIN3A}.

Finally, we note that the \ads\ $a_0, a_1, a_2, \ldots$ used in the definition of \sdcurva\ is crucial.  Although in principle one could define analogs of \sdcurva\ by using different sequences $a_0, a_1, a_2, \ldots$ in Equation~\ref{eqCURVDEF1}, such curvature functions would not satisfy all our properties.  A look at the proofs of our three theorems above will show that whereas Theorem~\ref{thmMAIN3A} would still hold with any sequence $a_0, a_1, a_2, \ldots$ in Equation~\ref{eqCURVDEF1}, Theorem~\ref{thmMAIN1} holds only if the sequence satisfies the recursive definition of the \ads\ given in Lemma~\ref{lemADS1} (3), and Theorem~\ref{thmMAIN2} holds because $a_n = 0$ for all odd $n$ (Lemma~\ref{lemADS1} (1)).  The author finds it somewhat remarkable that the \ads\ actually satisfies both these needed properties.

\section{Proofs}
\label{secPRO}

We start with a proof of Lemma~\ref{lemADS1}.

\demoname{Proof of Lemma~\ref{lemADS1}}
\noindent (1) This follows immediately from the fact that $B_n = 0$ for all odd integers $n$ such that $n \ge 3$ (see \cite[Sections~6.5]{G-K-P}).
\spce

\noindent (2) The fact that $a_0 = 1$ follows from the fact that $B_2 = 1/6$. To verify Equation~\ref{eqADS1}, we will need the following three formulas involving the Bernoulli numbers and the Bernoulli polynomials (denoted $B_k(x)$), the first two from \cite[Sections~6.5]{G-K-P}, and the third from \cite[p.\ 805]{A-S}.  For each non-negative integer $n$, we have
\begin{gather}
\sum_{i=0}^n \binom ni B_i = B_n\label{eqBB1}\\
\sum_{i=0}^n \binom ni B_i h^{n-i} = B_n(h)\label{eqBB2}\\
B_n(\frac 12) = -(1 - 2^{1-n})B_n\label{eqBB3}.
\end{gather}
By substituting $h = 1/2$ into Equation~\ref{eqBB2}, and then using Equation~\ref{eqBB3} and doing some rearranging, we obtain
\begin{equation}\label{eqBB4}
\sum_{i=0}^n \binom ni B_i 2^i = B_n (2 - 2^n).
\end{equation}
Next, substituting the definition of the $a_n$, given in Equation~\ref{eqBERN1}, into the left hand side of Equation~\ref{eqADS1}, and then using Equations~\ref{eqBB1} and \ref{eqBB4}, the fact that $B_0 = 1$ and $B_1 = -1/2$, and a standard identity involving binomial coefficients, it is a straightforward computation to show that Equation~\ref{eqADS1} holds.  We leave the details to the reader.
\edemoname

We now prove our three theorems, starting with Theorem~\ref{thmMAIN1}.

\demoname{Proof of Theorem~\ref{thmMAIN1}}
We compute
\begin{align*}
\allowdisplaybreaks
\sum_{\tau^p \in K} &(-1)^p \ka n{\tau^p} =  \sum_{p=0}^{n-2} \sum_{\tau^p \in K} (-1)^p \ka n{\tau^p}
 \intertext{\hfill using the fact that $\ka n{\tau^p} = 0$ when $p = n-1$ or $p = n$}\\
 &= \sum_{p=0}^{n-2} \sum_{\tau^p \in K} (-1)^p \left\{a_p\kmnna n{\tau^p} + \frac {a_p}2 \sum_{i=p+1}^{n-2} (-1)^{i-p} \sum_{\eta^i \offace \tau^p} \kmnna n{\eta^i}\right\}
 \intertext{\hfill using Equation~\ref{eqCURVDEF2}}\\
 &= \sum_{p=0}^{n-2} \sum_{\tau^p \in K} (-1)^p a_p\kmnna n{\tau^p} + \sum_{p=0}^{n-2} \sum_{\tau^p \in K} \sum_{i=p+1}^{n-2}\sum_{\eta^i \offace \tau^p}  (-1)^{i} \frac {a_p}2\kmnna n{\eta^i}\\
 &= \sum_{p=0}^{n-2} \sum_{\tau^p \in K} (-1)^p a_p\kmnna n{\tau^p} + \sum_{p=0}^{n-2} \sum_{i=p+1}^{n-2}\sum_{\eta^i \in K} \sum_{\tau^p \faceof \eta^i}  (-1)^{i} \frac {a_p}2\kmnna n{\eta^i}\\
 &= \sum_{p=0}^{n-2} \sum_{\tau^p \in K} (-1)^p a_p\kmnna n{\tau^p} + \sum_{p=0}^{n-2} \sum_{i=p+1}^{n-2}\sum_{\eta^i \in K}  (-1)^{i} \frac {a_p}2\binom {i+1}{p+1}\kmnna n{\eta^i}\\
 &= \sum_{i=0}^{n-2} \sum_{\eta^i \in K} (-1)^i a_i\kmnna n{\eta^i} + \sum_{i=1}^{n-2} \sum_{\eta^i \in K} \sum_{p=0}^{i-1} \binom {i+1}{p+1} (-1)^{i} \frac {a_p}2\kmnna n{\eta^i}\\
 &= \sum_{\eta^0 \in K} (-1)^0 a_0\kmnna n{\eta^0} + \sum_{i=1}^{n-2} \sum_{\eta^i \in K} (-1)^i \left\{a_i + \sum_{p=0}^{i-1} \binom {i+1}{p+1} \frac {a_p}2\right\}\kmnna n{\eta^i}\\
 &= \sum_{\eta^0 \in K} (-1)^0 \,\,1 \cdot \kmnna n{\eta^0} + \sum_{i=0}^{n-2} \sum_{\eta^i \in K} (-1)^i \,\,1 \cdot \kmnna n{\eta^i}
 \intertext{\hfill using Lemma~\ref{lemADS1} (2)}\\\noalign{\break}
 &= \sum_{i=0}^{n-2}  (-1)^i \sum_{\eta^i \in K}  \left\{\simrankk n{\eta^i} - \sum_{\sigma^n \offace \eta^i} \aang(\eta^i, \sigma^n)\right\}
 \intertext{\hfill using Equation~\ref{eqMOD1}}\\
 &= \eucn(K)
 \intertext{\hfill using the computation given in \cite[p.\ 387]{BL4}.}
\end{align*}
\edemoname

We now turn to the proof of Theorem~\ref{thmMAIN2}, starting with two lemmas.

\lem\label{lemPSEU}
Let $K$ be an \simpn n, where $n \ge 2$.  Suppose that each \nd{(n-1)}simplex of $K$ is the face of precisely two \nd{n}simplices.  Let $\tau^p$ be a 
\nd{p}simplex of $K$ with $0 \le p \le n-2$.  Then
\enum
\item[(i)] $\displaystyle \simrankk n{\tau^p} = 1$;
\item[(ii)] $\displaystyle 2f_{n-p-2}(\linkpl(\tau^p, K)) = (n-p)f_{n-p-1}(\linkpl(\tau^p, K))$.
\eenum
\elem

\demo
(i).  We use the discussion in Section~\ref{secSSC}.  Recall that we are thinking of simplices as open (and hence disjoint).  Using Remark~\ref{remCSET} (3), we know that $\tau^p$ is entirely contained in one of the sets $\cset nr(K)$.  Suppose $r \ne 2$.  By the same remark, we know that $\cset nr(K)$ is an 
$(n-1)$-manifold, and it is the union of simplices of $K$.  Therefore it must contain at least one \nd{(n-1)}simplex $\eta^{n-1}$.  By the definition of the set $\cset nr(K)$, it follows that $\eta^{n-1}$ must be the face of precisely $r$ \nd{n}simplices of $K$.  Because $r \ne 2$, we have a contradiction to the hypothesis of the lemma.  Hence we conclude that $r = 2$, and hence $\simrankk n{\tau^p} = 1$, using the definition of $\simrankk n{\tau^p}$.
\spce

\noindent (ii).  Let $\zeta^{n-p-1}$ be an \nd{(n-p-1)}simplex of $\linkpl(\tau^p, K)$.  We know that $\zeta^{n-p-1}$ has $n - p$ \nd{(n-p-1)}faces, and these faces are all in $\linkpl(\tau^p, K)$.  We claim that every \nd{(n-p-2)}simplex of $\linkpl(\tau^p, K)$ is the face of precisely two \nd{(n-p-1)}simplices of $\linkpl(\tau^p, K)$.  The desired result follows immediately from this claim.  To prove the claim, let $\eta^{n-p-2}$ be an \nd{(n-p-2)}simplex of $\linkpl(\tau^p, K)$.  Then $\tau^p \ast \eta^{n-p-2}$ is an \nd{(n - 1)}simplex of $K$.  It is seen that for every \nd{(n-p-1)}simplex $\alpha^{n-p-1}$ of $\linkpl(\tau^p, K)$ that has $\eta^{n-p-2}$ as a face, we obtain the \nd{n}simplex $\tau^p \ast \alpha^{n-p-1}$ of $K$ that has $\tau^p \ast \eta^{n-p-2}$ as a face.  Conversely, for every \nd{n}simplex $\beta^{n}$ of $K$ that has $\tau^p \ast \eta^{n-p-2}$ as a face, we can write $\beta^n$ as $\tau^p \ast \gamma^{n-p-1}$, where $\gamma^{n-p-1}$ is an \nd{(n-p-1)}simplex of $\linkpl(\tau^p, K)$ that has $\eta^{n-p-2}$ as a face.  Hence the number of \nd{(n-p-1)}simplices of $\linkpl(\tau^p, K)$ that have $\eta^{n-p-2}$ as a face equals the number of \nd{n}simplices of $K$ that have $\tau^p \ast \eta^{n-p-2}$ as a face.  By hypothesis the latter number is two, and hence there are also two \nd{(n-p-1)}simplices of $\linkpl(\tau^p, K)$ that have $\eta^{n-p-2}$ as a face.
\edemo

\lem\label{lemSIMP}
Let $n$ be an odd integer such that $n \ge 3$.  Let $\sigma^n$ be an \nd{n}simplex in Euclidean space, and let $\tau^p$ be a \nd{p}face of $\sigma^n$, where $p$ is an even integer such that $0 \le p \le n-2$.  Then
\begin{equation}\label{eqSIMP1}
\aang(\tau^p, \sigma^n) - \frac 12 \sum_{i=p+1}^{n-2} (-1)^{i+1} \!\!\!\sum_{\sigma^n \offace \eta^i \offace \tau^p} \aang(\eta^i, \sigma^n) = \frac 12 - \frac {n-p}4,
\end{equation}
where the inner summation is over all \nd{i}simplices $\eta^i$ of $\sigma^n$ that have $\tau^p$ as a face.
\elem

\demo
First, we show that Equation~\ref{eqSIMP1} is equivalent to the equation
\begin{equation}\label{eqSIMP2}
\sum_{i=p+1}^n (-1)^{i-p+1} \!\!\!\sum_{\sigma^n \offace \eta^i \offace \tau^p} \aang(\eta^i, \sigma^n) = 2\aang(\tau^p, \sigma^n).
\end{equation}
To see the equivalence, observe that by definition we have $\aang(\sigma^n, \sigma^n) = 1$, and that $\aang(\eta^{n-1}, \sigma^n) = 1/2$ for every \nd{(n-1)}face $\eta^{n-1}$ of $\sigma^n$.  It can be verified that there are precisely $n-p$ \nd{(n-1)}faces of $\sigma^n$ that contain $\tau^p$.
Then, using the fact that $n$ is odd and $p$ is even, we have 
\begin{equation}\label{eqSIMP3}
(-1)^{(n-1)-p+1} \sum_{\sigma^n \offace \eta^{n-1} \offace \tau^p} \aang(\eta^i, \sigma^n) = -\frac {n-p}2,
\end{equation}
and
\begin{equation}\label{eqSIMP4}
(-1)^{n-p+1} \sum_{\sigma^n \offace \eta^{n} \offace \tau^p} \aang(\eta^n, \sigma^n) = 1.
\end{equation}
Substituting Equations~\ref{eqSIMP3} and \ref{eqSIMP4} into Equation~\ref{eqSIMP2}, and doing some rearranging, yields Equation~\ref{eqSIMP1}.

Next, we need to verify that Equation~\ref{eqSIMP2} holds.  This equation follows from standard results concerning angle sums in convex polytopes and convex cones, and we sketch the proof.  We cite \cite{MC3} for some basic notation and results, though these results are found in many other sources as well; see \cite{GR1} for more background about convex polytopes.  

Withough loss of generality, we can translate $\sigma^n$ so that some point in $\tau^p$ is taken to the origin (recall that we are thinking of simplices as open).  Because all the angles under consideration are in $\sigma^n$, we can restrict our attention to the the affine span of $\sigma^n$, which we identify with $\rrr{n}$.

Let $Q$ be the polyhedral cone with apex at the origin generated by $\sigma^n$.  The set of apices of the cone $Q$, denoted $T$, is the affine span of $\tau^p$.  Clearly $T$ is a \nd{p}dimensional linear subspace $T$ of $\rrr{n}$  The faces of $Q$ correspond precisely to those faces of $\sigma^n$ that have $\tau^p$ as a face; for each face $\eta^i$ of $\sigma^n$ that has $\tau^p$ as a face, we denote its corresponding face in $Q$ by ${\hat \eta}^i$.  Moreover, we have $\aang({\hat \eta}^i, Q) = \aang(\eta^i, \sigma^n)$ for every appropriate $\eta^i$.  

Next, let $P = Q \cap T^\perp$, where $T^\perp$ is the \nd{(n-p)}dimensional linear subspace of $\rrr{n}$ that is perpendicular to $T$.  Then $P$ is an \nd{(n-p)}dimensional polyhedral cone, with one apex, namely the origin.  Moreover, it is seen that $Q = P \times T^\perp$.  The faces of $P$ correspond to the faces of $Q$; more precisely, for each \nd{i}face ${\hat \eta}^i$ of $Q$, there is a corresponding \nd{(i-p)}face ${\bar \eta}^i$ of $P$, where ${\bar \eta}^i = {\hat \eta}^i \cap T^\perp$ and ${\hat \eta}^i = {\bar \eta}^i \times T^\perp$.  Then, using Lemma~2 from \cite{MC}, and the fact that our angles are normalized, it is seen that for each face $\eta^i$ of $\sigma^n$ that has $\tau^p$ as a face, we have $\aang({\bar \eta}^i, P) = \aang({\hat \eta}^i, Q)$; hence $\aang({\bar \eta}^i, P) = \aang(\eta^i, \sigma^n)$.  Therefore, making use of the fact that $p$ is even, we see that in order to show Equation~\ref{eqSIMP2}, it suffices to show
\begin{equation}\label{eqSIMP5}
\sum_{i=1}^{n-p} (-1)^{i+1} \!\!\!\sum_{{\bar \eta}^i \faceof P} \aang({\bar \eta}^i, P) = 2\aang(0, P).
\end{equation}
However, Equation~\ref{eqSIMP5} is precisely the restatement for polyhedral cones of Sommerville's Theorem concerning the angle sums and volume of convex spherical polytopes, and hence Equation~\ref{eqSIMP5} is true.  See \cite[p.\ 157]{SOMM2} for the original statement of Sommerville's Theorem, and see \cite[p.\ 174]{MC3} for the restatement of this theorem for polyhedral cones (note that $n-p$ is odd, which is needed to deduce Equation~\ref{eqSIMP5} from Sommerville's Theorem).
\edemo

We are now ready for the proof of Theorem~\ref{thmMAIN2}.

\demoname{Proof of Theorem~\ref{thmMAIN2}}
Let $\tau^p$ be a \nd{p}simplex of $K$.  If $p = n - 1$ or $p = n$, then $\ka n{\tau^p} = 0$ by Remark~\ref{remCURV3}.  Now assume that $0 \le p \le n-2$.  If $p$ is odd, then $\ka n{\tau^p} = 0$ by Lemma~\ref{lemADS1} (1) combined with the definition of $\ka n{\tau^p}$.  From now on assume that $p$ is even.  We observe that the hypothesis on the links of even-dimensional simplices implies that every \nd{(n-1)}simplex of $K$ is the face of precisely two \nd{n}simplices.  Hence we can apply Lemma~\ref{lemPSEU} to $K$.

We compute
\begin{align*}
\allowdisplaybreaks
\frac 1{a_p} &\ka n{\tau^p} =  \kmnna n{\tau^p} + \frac 12 \sum_{i=p+1}^{n-2} (-1)^{i-p} \sum_{\eta^i \offace \tau^p} \kmnna n{\eta^i}\\
 &= \left[\simrankk n{\tau^p} - \sum_{\sigma^n \offace \tau^p} \aang(\tau^p, \sigma^n)\right] \\
&\qquad\qquad+ \frac 12 \sum_{i=p+1}^{n-2} (-1)^{i} \sum_{\eta^i \offace \tau^p} \left[\simrankk n{\eta^i} - \sum_{\sigma^n \offace \eta^i} \aang(\eta^i, \sigma^n)\right]\\
 &= 1 + \frac 12 \sum_{i=p+1}^{n-2} (-1)^{i} \sum_{\eta^i \offace \tau^p} 1\\
 &\qquad- \sum_{\sigma^n \offace \tau^p} \aang(\tau^p, \sigma^n) - \frac 12 \sum_{i=p+1}^{n-2} (-1)^{i} \sum_{\eta^i \offace \tau^p} \sum_{\sigma^n \offace \eta^i} \aang(\eta^i, \sigma^n)
\intertext{\hfill using Lemma~\ref{lemPSEU} (i), and doing some rearranging}\\\noalign{\break}
&= 1 - \frac 12 \sum_{k=0}^{n-p-3} \sum_{\omega^k \in \linkpl(\tau^p, K)} (-1)^{k}\\
 &\qquad- \sum_{\sigma^n \offace \tau^p} \aang(\tau^p, \sigma^n) - \frac 12 \sum_{\sigma^n \offace \tau^p} \sum_{i=p+1}^{n-2} (-1)^{i} \!\!\!\sum_{\sigma^n \offace \eta^i \offace \tau^p} \aang(\eta^i, \sigma^n)\\
&= 1 - \frac 12 \sum_{k=0}^{n-p-1} \sum_{\omega^k \in \linkpl(\tau^p, K)} (-1)^{k}  + \frac 12 \sum_{\omega^{n-p-2} \in \linkpl(\tau^p, K)} (-1)^{n-p-2}\\
 &\qquad\qquad\qquad\qquad+ \frac 12 \sum_{\omega^{n-p-1} \in \linkpl(\tau^p, K)} (-1)^{n-p-1}\\
 &\qquad- \sum_{\sigma^n \offace \tau^p} \left\{\aang(\tau^p, \sigma^n) - \frac 12 \sum_{i=p+1}^{n-2} (-1)^{i+1} \!\!\!\sum_{\sigma^n \offace \eta^i \offace \tau^p} \aang(\eta^i, \sigma^n)\right\}\\
 &= 1 - \frac 12 \chi(\linkpl(\tau^p, K))  - \frac 12 f_{n-p-2}(\linkpl(\tau^p, K))\\
 &\qquad\qquad+ \frac 12 f_{n-p-1}(\linkpl(\tau^p, K)) - \sum_{\sigma^n \offace \tau^p} \left\{\frac 12 - \frac {n-p}4\right\}
\intertext{\hfill using the fact that $n-p-2$ is odd, and Equation~\ref{eqSIMP1}}\\
 &= 1 - \frac 12 \cdot 2  - \frac 12 \frac {n-p}2 f_{n-p-1}(\linkpl(\tau^p, K))\\
 &\qquad\qquad+ \frac 12 f_{n-p-1}(\linkpl(\tau^p, K)) - f_{n-p-1}(\linkpl(\tau^p, K)) \left(\frac 12 - \frac {n-p}4\right)
\intertext{\hfill using the hypothesis of the theorem, and Lemma~\ref{lemPSEU} (ii)}\\
  &= 0,
\end{align*}
where the equality before the last uses the fact that $f_{n-p-1}(\linkpl(\tau^p, K))$ equals the number of \nd{n}simplices of $K$ that have $\tau^p$ as a face.
\edemoname

Finally, we have the proof of Theorem~\ref{thmMAIN3A}.

\demoname{Proof of Theorem~\ref{thmMAIN3A}}
To prove Equation~\ref{eqMAIN3A},  we first observe that if either $s$ or $p$ is odd, then it follows from Remark~\ref{remCURV3} (2) that both sides of Equation~\ref{eqMAIN3A} are zero, and so the equation holds.  From now on, assume that both $s$ and $p$ are even.  Clearly $s \le p$.

We start with the following preliminary.  Suppose that $\alpha^r$ is an \nd{r}simplex of $K$ that has $\zeta^p$ as a face, where $p+1 \le r \le n-2$.   We then claim that
\begin{equation}\label{eqQUEST2A}
\sum_{i=s+1}^{r} \sum_{\substack {\eta^i \in L \\ \eta^i \subseteq \alpha^r \\ \eta^i \offace \tau^s}}  (-1)^{i-s} = (-1)^{r-s}.
\end{equation}
To see why Equation~\ref{eqQUEST2A} holds, we first observe that 
\begin{equation}\label{eqQUEST3A}
\bigcup_{\substack {\eta^i \in L \\ \eta^i \subseteq \alpha^r \\ \eta^i \offace \tau^s}} \eta^i = |\starpl(\tau^s, L)| \cap \alpha^r,
\end{equation}
where, as always, we are thinking of simplices as open.  Let $T = \{\omega^k \in \linkpl(\tau^s, L) \mid \omega^k \subseteq \alpha^r\}$.  Observe that $|T|$ is an open \nd{(r-s-1)}disk.  We now have
\begin{equation}\label{eqQUEST4A}
\sum_{i=s+1}^{r} \sum_{\substack {\eta^i \in L \\ \eta^i \subseteq \alpha^r \\ \eta^i \offace \tau^s}}  (-1)^{i-s-1} = \sum_{k=0}^{r-s-1} \sum_{\omega^k \in T}  (-1)^{k}  = \sum_{k=0}^{r-s-1} (-1)^{k} f_k(T).
\end{equation}
We note that the sum $\sum_{k=0}^{r-s-1} (-1)^{k} f_k(T)$ is not necessarily equal to $\chi(T)$, because $T$ is not a simplicial complex.  However, we note that $\overline T$, the closure of $T$, is a simplicial complex (in particular, it is an \nd{(r-s-1)}disk).  Hence, using the discussion in \cite[Section~2]{BL4}, it is seen that 
\begin{equation}\label{eqQUEST5A}
\sum_{k=0}^{r-s-1} (-1)^{k} f_k(T) = (-1)^{r-s-1}.
\end{equation}
If we combine Equations~\ref{eqQUEST4A} and \ref{eqQUEST5A}, and multiply through by $-1$, we deduce Equation~\ref{eqQUEST2A}. 

One more preliminary observation.  Suppose that $\nu^i$ is an \nd{i}simplex of $L$ that is contained in a (unique) \nd{r}simplex $\mu^r$ of $K$.  We observe that $\simrankk n{\nu^i} = \simrankk n{\mu^r}$ and $\kmnna n{\nu^i} = \kmnna n{\mu^r}$, where the left hand side of each of these equations is with respect to $L$, and the right hand side is with respect to $K$ (in general, if we write ``$\kmnna n{\eta^i}$,'' that is with respect to the simplicial complex of which $\eta^i$ is a simplex, which is always clear from context).  

We can now demonstrate Equation~\ref{eqMAIN3A}, by computing
\begin{align*}\label{eqINVAR1}
\allowdisplaybreaks
a_p\ka n{\tau^s} &= a_pa_s\kmnna n{\tau^s} + a_p\frac {a_s}2 \sum_{i=s+1}^{n-2} (-1)^{i-s} \sum_{\substack {\eta^i \in L \\ \eta^i \offace \tau^s}} \kmnna n{\eta^i}\\
 &= a_pa_s\kmnna n{\zeta^p} + \frac {a_pa_s}2 \sum_{\substack {\eta^i \in L \\ s+1 \le i \le n-2 \\ \eta^i \offace \tau^s}} (-1)^{i-s} \kmnna n{\eta^i}\\\noalign{\break}
 &= a_pa_s\kmnna n{\zeta^p} + \frac {a_pa_s}2 \sum_{\substack {\alpha^r \in K \\ p+1 \le r \le n \\ \alpha^r \offace \zeta^p}} \sum_{\substack {\eta^i \in L \\ \eta^i \subseteq \alpha^r \\ s+1 \le i \le \min\{r, n-2\} \\ \eta^i \offace \tau^s}} (-1)^{i-s} \kmnna n{\alpha^r}\\
 &= a_pa_s\kmnna n{\zeta^p} + \frac {a_pa_s}2 \sum_{r=p+1}^{n} \sum_{\substack {\alpha^r \in K \\ \alpha^r \offace \zeta^p}} \kmnna n{\alpha^r} \sum_{i=s+1}^{r} \sum_{\substack {\eta^i \in L \\ \eta^i \subseteq \alpha^r \\ \eta^i \offace \tau^s}}  (-1)^{i-s}\\
\intertext{\hfill because $\kmnna n{\alpha^r} = 0$ for $r = n-1$ or $r = n$}\\
 &= a_sa_p\kmnna n{\zeta^p} + a_s\frac {a_p}2 \sum_{r=p+1}^{n-2} \sum_{\substack {\alpha^r \in K \\ \alpha^r \offace \zeta^p}} \kmnna n{\alpha^r} (-1)^{r-p} 
\intertext{\hfill Using Equation~\ref{eqQUEST2A}, and the fact that both $s$ and $p$ are even}\\
 &= a_s\ka n{\zeta^p}.
\end{align*}

Finally, to deduce Equation~\ref{eqMAIN3B}, assume $s = p$.  There are two cases to consider.  When $a_s = a_p = 0$, then both sides of Equation~\ref{eqMAIN3B} are zero by definition.  When $a_s = a_p \ne 0$, then Equation~\ref{eqMAIN3B} follows from Equation~\ref{eqMAIN3A}.
\edemoname

\providecommand{\bysame}{\leavevmode\hbox to3em{\hrulefill}\thinspace}


\begin{thebibliography}{McM86}

\bibitem[AS72]{A-S}
M.~Abramowitz and I.~Stegun, \emph{Handbook of mathematical functions with
  formulas, graphs, and mathematical tables}, Applied Mathematics Series,
  vol.~55, National Bureau of Standards, Washington, DC, 1972.

\bibitem[Ban67]{BA1}
Thomas Banchoff, \emph{Critical points and curvature for embedded polyhedra},
  J. Diff. Geom. \textbf{1} (1967), 245--256.

\bibitem[Ban70]{BA2}
Thomas Banchoff, \emph{Critical points and curvature for embedded polyhedral
  surfaces}, Amer. Math. Monthly \textbf{77} (1970), 475--485.

\bibitem[Ban83]{BA3}
Thomas Banchoff, \emph{Critical points and curvature for embedded polyhedra,
  ii}, Progress in Math. \textbf{32} (1983), 34--55.

\bibitem[Blo]{BL10}
Ethan~D. Bloch, \emph{Critical points and the angle defect}, to appear.

\bibitem[Blo98]{BL4}
Ethan~D. Bloch, \emph{The angle defect for arbitrary polyhedra}, Beitr{\"a}ge
  Algebra Geom. \textbf{39} (1998), 379--393.

\bibitem[Bud89]{BUDA}
Lothar Budach, \emph{{L}ipschitz-{K}illing curvatures of angular partially
  ordered sets}, Adv. in Math. \textbf{78} (1989), 140--167.

\bibitem[Che83]{CHEE}
J.~Cheeger, \emph{Spectral geometry of singular riemannian spaces}, J. Diff.
  Geom. \textbf{18} (1983), 575--657.

\bibitem[CMS84]{C-M-S}
J.~Cheeger, W.~Muller, and R.~Schrader, \emph{On the curvature of piecewise
  flat spaces}, Commun. Math. Phys. \textbf{92} (1984), 405--454.

\bibitem[Fed82]{FE}
P.~J. Federico, \emph{Descartes on polyhedra}, Springer-Verlag, New York, 1982.

\bibitem[GKP94]{G-K-P}
Ronald Graham, Donald Knuth, and Oren Patashnik, \emph{Concrete mathematics},
  second ed., Addison-Wesley, Reading, MA, 1994.

\bibitem[Gr{\"u}67]{GR1}
Branko Gr{\"u}nbaum, \emph{Convex polytopes}, John Wiley \& Sons, New York,
  1967.

\bibitem[Gr{\"u}68]{GR2}
Branko Gr{\"u}nbaum, \emph{{G}rassman angles of convex polytopes}, Acta. Math.
  \textbf{121} (1968), 293--302.

\bibitem[GS91]{G-S2}
Branko Gr{\"u}nbaum and G.~C. Shephard, \emph{{D}escartes' theorem in $n$
  dimensions}, Enseign. Math. (2) \textbf{37} (1991), 11--15.

\bibitem[Gui59]{GUIN}
A.~P. Guinand, \emph{A note on the angles in an $n$-dimensional simplex}, Proc.
  Glasgow Math. Assoc. \textbf{4} (1959), 58--61.

\bibitem[Hud69]{HU}
J.~F.~P. Hudson, \emph{Piecewise linear topology}, Benjamin, Menlo Park, CA,
  1969.

\bibitem[Kno28]{KNOP}
Konrad Knopp, \emph{Theory and application of infinite series}, Blackie and
  Son, London, 1928.

\bibitem[McM75]{MC}
Peter McMullen, \emph{Non-linear angle-sum relations for polyhedral cones and
  polytopes}, Math. Proc. Cambridge Philos. Soc. \textbf{78} (1975), 247--261.

\bibitem[McM86]{MC3}
Peter McMullen, \emph{Angle-sum relations for polyhedral sets}, Mathematika
  \textbf{33} (1986), 173--188.

\bibitem[MS74]{M-S}
John~W. Milnor and James Stasheff, \emph{Characteristic classes}, Ann. of Math.
  Studies, vol.~76, Princeton U. Press, Princeton, NJ, 1974.

\bibitem[Mun84]{MU3}
James~R. Munkres, \emph{Elements of algebraic topology}, Addison-Wesley, Menlo
  Park, CA, 1984.

\bibitem[Pes56]{PESC}
E.~Peschl, \emph{{W}inkelrelationen am {S}implex und die {E}ulersche
  {C}harakteristik}, Bayer. Akad. Wiss., Math.-Nat. Kl.S.-B. \textbf{1955}
  (1956), 319--345.

\bibitem[PS67]{P-S}
M.~A. Perles and G.~C. Shephard, \emph{Angle sums of convex polytopes}, Math.
  Scand. \textbf{21} (1967), 199--218.

\bibitem[She67]{SH}
G.~C. Shephard, \emph{An elementary proof of {G}ram's theorem for convex
  polytopes}, Canad. J. Math. \textbf{19} (1967), 1214--1217.

\bibitem[She68]{SH2}
G.~C. Shephard, \emph{Angle deficiencies of convex polytopes}, J. London Math.
  Soc. \textbf{43} (1968), 325--336.

\bibitem[Som]{SOMM2}
D.~M.~Y. Sommerville, \emph{An introduction to the geometry of $n$ dimensions},
  Methuen, London.

\bibitem[Win82]{WINT}
P.~Wintgen, \emph{Normal cycle and integral curvature for polyhedra in
  riemannian manifolds}, Differential Geometry (Amsterdam) (Gy. So\'os and
  J.~Szenthe, eds.), vol.~31, North-Holland, Amsterdam, 1982, pp.~805--816.

\bibitem[Z{\"a}h]{ZA}
M.~Z{\"a}hle, \emph{Approximation and characterization of generalized
  {L}ipschitz-{K}illing curvatures}, Ann. Global Annal. Geom. \textbf{8} (1990), 249--260.

\end{thebibliography}
\end{document}